\documentclass[doublespacing]{elsart}
\usepackage{amsmath}
\usepackage{amsfonts}   
\usepackage{amssymb}

\journal{}

\newcommand{\mlabel}[1]{\label{#1}
}

\newcommand{\seq}{\begin{equation}}                 
\newcommand{\eeq}[1]{\label{#1}\end{equation}
    }

\newcommand{\epf}{$ \quad \Box$ \par \vspace{1ex}}

\newtheorem{Theorem}{Theorem}[section]
\newcommand{\sthm}{\begin{Theorem}}         
\newcommand{\ethm}{\end{Theorem}}           

\newtheorem{Corollary}[Theorem]{Corollary}
\newcommand{\scor}{\begin{Corollary}}       
\newcommand{\ecor}{\end{Corollary}}         

\newtheorem{Lemma}[Theorem]{Lemma}
\newcommand{\slm}{\begin{Lemma}}            
\newcommand{\elm}{\end{Lemma}}              

\newtheorem{Example}[Theorem]{\sc Example}
\newcommand{\sex}{\begin{Example}\rm}        
\newcommand{\eex}{\end{Example}}             

\newcommand{\seql}{\begin{eqnarray*}}       
\newcommand{\eeql}{\end{eqnarray*}}

\newcommand{\smlist}[1]{\begin{list}           
                      {(#1{zzcount})}{\usecounter{zzcount}}}
\newcommand{\elist}{\end{list}}

\newcommand{\for}{\quad \mbox{for} \quad}

\newcommand{\G}{\Gamma}
\renewcommand{\l}{\lambda}

\newcommand{\s}{\sigma}

\renewcommand{\O}{\Omega}

\providecommand{\norm}[1]{\lVert#1\rVert}
\providecommand{\abs}[1]{\lvert#1\rvert}

\def\R{{\Bbb R}}

\def\e{\varepsilon}

\def\lb{\label}

\begin{document}
\begin{frontmatter}

\title{Positive periodic solutions of singular systems}
\author{Haiyan Wang}
\ead{wangh@asu.edu}
\address{Division of Mathematical and Natural Sciences\\ Arizona State University \\ Phoenix, AZ 85069-7100, USA}

\pagenumbering{arabic}

\begin{abstract}The existence and multiplicity of positive
periodic solutions for second order non-autonomous singular dynamical systems are established
with superlinearity or sublinearity assumptions at infinity for an appropriately chosen parameter.
Our results provide a unified treatment for the problem and significantly improve several results in
the literature. The proof of our results is based on the Krasnoselskii fixed point theorem in a cone.
\end{abstract}


\end{frontmatter}

\section{Introduction}

\quad\quad In a recent series of papers, Chu and Torres \cite{ct}, Chu, Torres and Zhang \cite{chujde2007}, Franco and Webb \cite{fw},
Franco and Torres\cite{franzoTo}, Jiang, Chu and  Zhang \cite{jcz}, Torres \cite{t03,t04,t06}, the existence and multiplicity of
positive periodic solutions for the singular systems
\begin{equation}\label{eq1}
\ddot{x}+a(t)x= f(t,x)+ e(t)
\end{equation}
and
\begin{equation}\label{eq2}
-\ddot{x}+a(t)x= f(t,x)+ e(t)
\end{equation}
have been studied,
where $a(t),e(t) \in C(\R,\R^{n})$, $f(t,x)\in C(\R\times(\R^{n}\backslash\{0\}),\R^{n})$ are T-periodic in $t$ with a singularity at $x=0$,
$$
\lim_{ x \to 0} f_i(t,x)=\infty, i=1,...,n.
$$
(\ref{eq1}) and (\ref{eq2}) represent singularities
of repulsive type and attractive type respectively. One closely
related example of the above systems is
\begin{equation}\label{c}
\ddot{x}+a x+ \nabla_{x}V(t,x)=e(t)
\end{equation}
with $V(t,x)=(\frac{1}{\sqrt{\sum x_i^2}})^{\alpha+1}, \alpha>0$, which was studied in \cite{Majer1991}. A positive periodic solution of the above systems is of interest because it is a non-collision
periodic orbit of the singular systems. Periodic solutions of singular systems has been studied over many years, see,
for example, \cite{az, bc, ct, chujde2007, dmm, fer, fw, franzoTo, g, jcoa, jcz, LS} and \cite{Majer1991,rtv,rt, s90, te, t03, t04, t06, z99}. One of the common assumptions
to guarantee the existence of  is a so-called strong force assumption ( corresponds to the case $\alpha \geq 1$ in (\ref{c})), see, for example, \cite{az,g} and references therein.  However, more recently,
the existence of positive periodic solutions of the singular systems has been established with a weak force
condition \cite{ct, chujde2007, fw,  franzoTo, rtv, rt, t04, t06}.

The variational arguments have been the most used techniques to deal with the problem, see, for example, \cite{az,Majer1991,s90,ta, te}. More recently,
the method of lower and upper solutions, the Schauder's fixed
point theorem and the Krasnoselskii fixed point theorem in a cone have been employed to investigate the existence
of positive periodic solutions of the systems \cite{bc, ct, chujde2007, fw, franzoTo,jcoa, jcz,OREGANWANG2006,t03,t04,t06}. There is a rich literature on the use of the Krasnoselskii fixed point theorem for
the existence of positive solutions of boundary value problems for general second-order differential equations (refer to \cite{EW94,EHW,HW} and many other  papers).

Motivated by these recent developments, we investigate the existence and multiplicity of positive periodic solutions of the singular systems
by the Krasnoselskii fixed point theorem.
In this paper, we are able to obtain several existence results
based on the Krasnoselskii fixed point theorem by constructing
a cone defined on a product space.  Similar cones have been proposed to study the existence of
positive solutions of boundary value problems for systems of differential equations in several papers of the author and
his co-authors \cite{DUNWANG2,DUNWANG1,HWJMAA1}. We also note a related cone is used to
study the existence of positive periodic solutions of singular periodic systems \cite{fw,t04}.
It seems that the Krasnoselskii fixed point theorem on compression and expansion of cones is quite effective in dealing with the problem.
In fact, by choosing appropriate cones, the singularity of the systems is essentially removed and the associated operator becomes well-defined for
certain ranges of functions even when $e_i$ is negative.

%
%
%
%

This paper is organized  as follows. Main results are given in Section 2.  In Section 3, we define a cone
and discuss several properties of the equivalent operator on the cone. In order to simplify the proof in Section 3, we establish a series of lemmas and corollaries to estimate the operator.
All the corollaries are the corresponding results for $e_i$ taking negative values. The proof of the main results is presented in Sections 4 and 5.

\section{Main results}

In this section, we present our main results for  the existence and multiplicity of positive periodic solutions of singular systems of repulsive type (\ref{eq1}).
For (\ref{eq2}), all the results can be proved in the  same way. First, we state a condition to guarantee the positiveness of the Green's function of the following scalar problems, $i=1, 2,\cdots, n,$
\begin{equation}\label{p1}
  x_i'' +a_i(t) x_i=e_i(t)
\end{equation}
with periodic boundary conditions
$ x_i(0)=x_i(T), \quad x'_i(0)=x'_i(T),$
where $x=(x_1,x_2,\cdots,x_n),$  and $a_1,a_2,\cdots,a_n$ and $e_1,e_2,\cdots,e_n$ are $T$-periodic continuous functions.
Let $G_i(t,s) \in C([0, T], \R) $ be the Green functions associated with (\ref{p1}). Now the periodic solution
$x(t)=(x_1(t),x_2(t),\cdots,x_n(t))$ of $(\ref{p1})$ is given by
$$x_i(t)= \int_{0}^{T}G_i(t,s)e_i(s) ds.$$
When $a_i(t)=k^2,$ $0 < k< \frac{\pi}{T}$,  the Green function $G_i$ takes the following form,
\[G_i(t,s)=\left\{\begin{array}{l}
\frac{\sin k(t-s)+\sin k(T-t+s)}{2k(1-\cos kT)},\quad 0\leq s\leq t\leq T,\\
\frac{\sin k(s-t)+\sin k(T-s+t)}{2k(1-\cos kT)},\quad 0\leq
t\leq s\leq T.\end{array}\right.\]
We can verify that $G_i$ is strictly positive. In fact, let $\hat{G}(x)=\frac{\sin (kx)+\sin k( T-x)}{2k(1-\cos k T)},
x \in [0,  T]$. It is easy to check that $\hat{G}$ is increasing
on $[0, \frac{T}{2}]$ and decreasing on $[\frac{T}{2},  T]$, and
$G(t,s)=\hat{G}(|t-s|)$. Thus
$$
0< \frac{\sin k T}{2k(1-\cos k T)}= \hat{G}(0) \leq G(t,s) \leq
\hat{G}(\frac{T}{2})= \frac{\sin\frac{kT}{2}}{k(1-\cos k T)}=\frac{1}{2k\sin \frac{kT}{2}}
$$
for $ s, t \in [0,  T ].$ The same estimates can also be found in \cite{fw,OREGANWANG2006,t03}.  For a non-constant function $a_i(t)$, there is a criterion
discussed in \cite{t03,zl} to guarantee the positiveness of the Green's functions. Therefore, we always assume the following assumption (A) is true
for systems of repulsive type (\ref{eq1}) throughout the paper.
\begin{itemize}
\item[(A)] The Green function $G_i(t,s),$ associated with (\ref{p1}), is positive for
all $(t,s)\in [0,T]\times [0,T]$, $i=1, 2,\cdots, n.$
\end{itemize}
Under hypothesis (A), we denote
\begin{equation}\lb{d}
\begin{split}
0< m_i&=\min_{0\leq s,t\leq T}G_i(t,s), \, M_i=\max_{0\leq s,t\leq T}G_i(t,s),\\
 0<\sigma_i &=\frac{m_i}{M_i}, \, \s=\min_{i=1,...,n}\{\s_i\}>0.
\end{split}
\end{equation}

We now examine the existence and multiplicity of positive periodic solutions of the following form, for $i=1,...,n$
\begin{equation}\label{eq3}
\ddot{x_i}+a_i(t)x_i=\l g_i(t) f_i(x)+\l e_i(t).
\end{equation}
with $\l>0$ is a positive parameter.  By a positive $T$-periodic solution, we
mean a positive $T$-periodic function in $C^2(\R,\R^{n})$ solving corresponding systems and each component is positive for all $t$. Let $\mathbb{R}_+=[0, \infty)$, $\mathbb{R}_+^n=\Pi_{i=1}^n \mathbb{R}_+$, and denote by $\abs{x}=\sum_{i=1}^n\abs{x_i}$ the usual norm of $\mathbb{R}_+^n$ for
$x= (x_1,\cdots, x_n) \in \mathbb{R}^n$. We will make the following assumptions

\noindent
(H1)~~ $f_i(x)$ is a scalar continuous function defined for $\abs{x} > 0,$ and $f_i(x)>0$ for $\abs{x}>0$,  $i=1,\dots,n$. \\
(H2)~~ $a_i(t), g_i(t), e_i(t)$ are $T$-periodic continuous scalar functions in $t \in \R$, $a_i(t), g_i(t)
\geq 0, t \in [0, T], \int^T_0 g_i(t)dt>0$, $i=1,\dots,n.$

We state our first theorem as follows.

\sthm\mlabel{th3} Let \rm{(A),(H1),(H2)} hold,  and $e_i(t) \geq 0$ for $t \in [0, T], i=1,...,n.$ Assume that $\lim_{ x \to 0} f_i(x)=\infty, i=1,...,n.$ \\
(a). If $\lim_{ \abs{x} \to \infty} \frac{f_i(x)}{\abs{x}}=0$,  $i=1,\dots,n$ , then, for all $\l>0$, (\ref{eq3}) has
a positive periodic solution. \\
(b). If $\lim_{ \abs{x} \to \infty} \frac{f_i(x)}{\abs{x}}=\infty$  for $i=1,\dots,n$, then, for all sufficiently small $\l>0$, (\ref{eq3}) has
two positive periodic solutions. \\
(c). There exists a $\l_0>0$ such that (\ref{eq3}) has a positive periodic solution for $0 < \l < \l_0$.
\ethm

When $e_i(t)$ takes negative values, we give the following theorem.  We need a stronger condition on $g_i$.\\
(H3)~~ $g_i(t)>0$ for $t \in [0, T]$, $i=1,\dots,n.$
\sthm\mlabel{th5} Let \rm{(A),(H1),(H2), (H3)} hold. Assume that $\lim_{ x \to 0} f_i(x)=\infty, i=1,...,n.$ \\
(a). If $\lim_{ \abs{x} \to \infty} f_i(x)=\infty$ and  $\lim_{ \abs{x} \to \infty} \frac{f_i(x)}{\abs{x}}=0$,  $i=1,\dots,n$ , then there exists  $\l_0>0$ such that
(\ref{eq3}) has a positive periodic solution for $\l > \l_0$. \\
(b). If $\lim_{ \abs{x} \to \infty} \frac{f_i(x)}{\abs{x}}=\infty$  for  $i=1,\dots,n$ , then, for all sufficiently small $\l>0$, (\ref{eq3})  has
two positive periodic solutions.\\
(c). There exists a $\l_1>0$ such that (\ref{eq3})  has a positive periodic solution for $0 < \l < \l_1$. \\
\ethm

Now we apply Theorems \ref{th3}, \ref{th5} to the following two-dimensional singular system,
which has been examined in \cite{chujde2007,franzoTo,jcoa}.
\begin{equation} \label{example}
\left\{ \begin{aligned}
\ddot{x}+a_1(t)x=\l \Big(\sqrt{x^2+y^2}\Big)^{-\alpha}+\l \Big(\sqrt{x^2+y^2}\Big)^{\beta}+\l e_1(t) \\
\ddot{y}+a_2(t)y=\l \Big(\sqrt{x^2+y^2}\Big)^{-\alpha}+\l \Big(\sqrt{x^2+y^2}\Big)^{\beta}+\l e_2(t),
\end{aligned} \right.
\end{equation}
with $\alpha, \beta>0$, $a_1 \geq 0,a_2 \geq 0, e_1, e_2$ are $T$-periodic continuous in $t$. We only need to note the following inequality
$$\sqrt{x^2+y^2} \leq \abs{x}+\abs{y} \leq \sqrt{2}\sqrt{x^2+y^2}$$
since we use the summation norm in our theorems.
For nonnegative $e_1,e_2$, Corollary \ref{cor1} is an application of Theorem \ref{th3}.
\scor\mlabel{cor1}. Assume that $a_1,a_2, e_1, e_2$ are $T$-periodic continuous in $t$ and that $a_1,a_2$ satisfy
the assumption (A). Also assume that $e_1\geq 0$ and $e_2\geq 0 $ for $t \in [0,T]$. Let $\alpha>0,\beta>0$, $\l>0$.\\
(a). If $0< \beta< 1$, then, for all $\l>0$, (\ref{example}) has
a positive periodic solution. \\
(b). If $\beta >1$, then, for all sufficiently small $\l>0$, (\ref{example})  has
two positive periodic solutions. \\
(c). There exists a $\l_0>0$ such that (\ref{example}) has a positive periodic solution for $0 < \l < \l_0$.
\ecor
When $e_1, e_2$ take negative values, we have the following corollary from Theorem \ref{th5}.
\scor\mlabel{cor2}. Assume that $a_1,a_2, e_1, e_2$ are $T$-periodic continuous in $t$, and that $a_1,a_2$
satisfy  the assumption (A). Let $\alpha>0,\beta>0$ and $\l>0$. \\
(a). If $0< \beta <1$ , then there exists  $\l_0>0$ such that
(\ref{example}) has a positive periodic solution for $\l > \l_0$. \\
(b). If $\beta>1$ , then, for all sufficiently small $\l>0$, (\ref{example})  has
two positive periodic solutions.\\
(c). There exists a $\l_1>0$ such that (\ref{example})  has a positive periodic solution for $0 < \l < \l_1$. \\
\ecor

We remark that the conclusions (b) of Theorems \ref{th3}, \ref{th5} are still valid if at least one component of
$f$ satisfies $\lim_{ \abs{x} \to \infty} \frac{f_i(x)}{\abs{x}}=\infty$. In addition,
analogous results are true if one considers a system that not every component is singular at zero.
For simplicity, every component of $f(t,x)$ is assumed to be singular
at zero in this paper.
Also we comment that Theorems \ref{th3} and \ref{th5} can be extended to the following more general system
\begin{equation}\label{eq6}
\ddot{x_i}+a_i(t)x_i=\l f_i(t,x)+\l e_i(t)
\end{equation}
if $f_i(t,x)$ satisfies $p_i(t) h_i(x) \leq f_i(t,x) \leq q_i(t)H_i(x), i=1, ..., n$ with appropriate conditions on $p_i,h_i, q_i,H_i$.

In comparison with some related results in \cite{ct, chujde2007, fw, franzoTo, jcz, t03, t04, t06},
the existence and multiplicity results in this paper can be applied to any periodic continuous function $e_i$. Of course, our results require the parameter $\l$ sufficiently small or large.
From Corollaries \ref{cor1} and \ref{cor2}, for $\alpha>0$, (\ref{example}) always has a positive periodic solution(s) if the parameter $\l$
is appropriately chosen according to $1>\beta>0$ or $\beta>1$. These results further suggest that both a strong force assumption and weak singularity
contribute to the existence of a positive solution(s) as long as certain conditions are met. Also it should be pointed out
that, for the non-singular case ($\alpha \leq 0$), several possible combinations of superlinear and sublinear assumptions at zero and infinity were considered in
\cite{OREGANWANG2006} to obtain one or two positive periodic solutions of periodic boundary value problems. Finally, we provide a unified treatment of the
problem for several important cases, and the conditions of our theorems are quite easy to verify.

We have formulated our arguments in a series of lemmas and corollaries to avoid repeated
arguments in the proofs of the results. All the corollaries in Section \ref{preli} are the corresponding results for $e_i$ which may take negative values.
It seems, to some extend, that the lemmas and corollaries themselves are of importance, and reveal significant properties
of the singular systems. We hope that they can be used in future research.

\section{Preliminary results}\label{preli}

We recall some concepts and conclusions of an operator in a cone. Let $E$
be a Banach space and $K$ be a closed, nonempty subset of $E$. $K$
is said to be a cone if $(i)$~$\alpha u+\beta v\in K $ for all
$u,v\in K$ and all $\alpha,\beta \geq 0$ and $(ii)$~$u,-u\in K$ imply
$u=0$. The following well-known result of the fixed
point theorem is crucial in our arguments.

\slm\mlabel{lm1} {\rm (\cite{KRAS})} Let $X$ be a Banach
space and $K\ (\subset X)$ be a cone. Assume that $\Omega_1,\
\Omega_2$ are open subsets of $X$ with $0 \in \Omega_1,\bar\Omega_1 \subset \Omega_2$, and let
$$
\mathcal{T}: K \cap (\bar{\Omega}_2\setminus \Omega_1 ) \rightarrow K
$$
be completely continuous operator such that either
\begin{itemize}
\item[{\rm (i)}] $\| \mathcal{T}u \| \geq \| u \|,\ u\in K\cap \partial
     \Omega_1$ and $ \| \mathcal{T}u \| \leq \| u \|,\ u\in K\cap \partial
     \Omega_2$; or

\item[{\rm (ii)}] $\| \mathcal{T}u \| \leq \| u \|,\ u\in K\cap \partial
     \Omega_1$ and $\| \mathcal{T}u \| \geq \| u \|,\ u\in K\cap \partial
     \Omega_2$.
\end{itemize}
Then $\mathcal{T}$ has a fixed point in $K \cap ( \bar \Omega_2 \backslash
     \Omega_1)$.

\elm
Consider the Banach space $X =\underbrace{C[0, T] \times ... \times C[0, T]}_{n}$, and for
$x=(x_1,...,x_n) \in X,$ let
$$\displaystyle{\norm{x}= \sum_{i=1}^n \sup_{t\in[0, T]} \abs{x_i(t)}}.$$
Denote by $K$ the cone
\begin{equation*}
\begin{split}
K = \{&x=(x_1,...,x_n) \in X: x_i(t)\geq 0,\; t \in [0, T], \; i=1,...,n,\\
&\rm{and} \; \min\limits_{0 \leq t \leq  T }\sum_{i=1}^nx_i(t) \geq \s \norm{x}\}
\end{split}
\end{equation*}
where $\s$ is defined in (\ref{d}). Also, for $r>0$, let
$$
\O_r  = \{x \in K: \norm{x} < r \}.
$$
Note that $\partial \O_r = \{x \in K: \norm{x}=r\}$.

Let us define $\mathcal{T}_{\l}=(\mathcal{T}_{\l}^1,...,\mathcal{T}_{\l}^n): K\setminus \{0\} \to X$, where $\mathcal{T}_{\l}^i$, $i=1,...,n$, are
\begin{equation}\label{T_def}
\mathcal{T}_{\l}^ix(t) =\l\int^{T}_0  G_i(t,s)\Big( g_i(s)f_i(x(s))+ e_i(s)\Big)ds , \;\; 0 \leq t \leq  T.
\end{equation}
When $e_i$ is nonnegative, $g_i(s)f_i(x(s))+ e_i(s)$ is nonnegative. If $e_i$ takes negative values, we will choose $x(s)$
so that  $g_i(s)f_i(x(s))+ e_i(s)$ is
nonnegative. This is possible because $\lim_{ x \to 0} f_i(x)=\infty$ or $\lim_{ \abs{x} \to \infty} f_i(x)=\infty$.

Now if $x$ is a fixed point of $\mathcal{T}_{\l}$ in  $K\setminus\{0\}$, then $x$ is a positive solution of (\ref{eq3}).
Also note that each component $x_i(t)$ of any nonnegative periodic solution $x$ is strictly positive for all $t$
because of the positiveness of the Green functions and assumptions (H1) and (H2). We now look at several properties of the operator.


\slm\mlabel{lm-compact}
Assume \rm{(A),(H1),(H2)} hold and $e_i(t)\geq 0, t \in [0,T], i=1,...,n$. Then  $\mathcal{T} _{\l}( K \setminus \{0\}) \subset K$ and $\mathcal{T}_{\l}:K \setminus \{0\} \to K$ is completely continuous.
\elm
\pf If $x \in K \setminus \{0\}$, then $\min_{ t \in [0,T]} \sum_{i=1}^n |x_i(t)| \geq \s \norm{x}>0$, and then $\mathcal{T} _{\l}$ is defined.
Now we have that, for  $i=1,\dots,n$
\begin{equation*}
\begin{split}
 \min_{ t \in [0,T]} \sum_{i=1}^n \mathcal{T}_{\lambda}^i x(t) & \geq \sum_{i=1}^n \min\limits_{0 \leq t \leq T}\mathcal{T}_{\l}^ix(t)\\
 & \geq \sum_{i=1}^n m_i \l \int^{ T}_{0}\big( g_i(s)f_i(x(s))+ e_i(s)\big)ds\\
 & = \sum_{i=1}^n \s_i \l M_i\int^{ T}_{0}\big( g_i(s)f_i(u(s))+ e_i(s)\big)ds  \\
 & \geq \sum_{i=1}^n \s_i \sup_{0 \leq t \leq T}{\mathcal{T}_{\l}^ix(t)}\\
 &  \geq \s \sum_{i=1}^n \sup_{0 \leq t \leq T}{\mathcal{T}_{\l}^ix(t)}=\s \norm{\mathcal{T}_{\l}x}.
\end{split}
\end{equation*}
Thus, $\mathcal{T} _{\l}( K \setminus \{0\} \subset K$.  It is easy to verify that $\mathcal{T} _{\l}$ is
completely continuous.
\epf
If $e_i$ takes negative values, we need to choose appropriate domains so that $g_i(s)f_i(x(s))+ e_i(s)$ becomes
nonnegative. The proof of $\mathcal{T} _{\l}( K \setminus \{0\}) \subset K$ and $\mathcal{T} _{\l}( K\setminus \O_R) \subset K$ in
Corollary \ref{cordefinition} is the same as in Lemma \ref{lm-compact}.
\scor\mlabel{cordefinition}
Assume \rm{(A),(H1),(H2), (H3)} hold. \\
a) If $\lim_{ x \to 0}  f_i(x)=\infty, i=1,\dots,n$, there is a $\delta>0$ such that if $0<r<\delta$, then   $\mathcal{T} _{\l}$ is defined on $\bar{\O}_r \setminus \{0\}$, $\mathcal{T} _{\l}(\bar{\O}_r \setminus \{0\}) \subset K$
 and $\mathcal{T}_{\l}:\bar{\O}_r \setminus \{0\} \to K$ is completely continuous.\\
b) If $\lim_{ x \to \infty}  f_i(x)=\infty, i=1,\dots,n$, there is a $\Delta>0$ such that if $R>\Delta$, then $\mathcal{T} _{\l}$ is defined on $K\setminus \O_R$,  $\mathcal{T} _{\l}( K\setminus \O_R) \subset K$ and $\mathcal{T}_{\l}:K\setminus \O_R \to K$ is completely continuous.
\ecor
\pf We split $g_i(s) f_i(x(s))+e_i(t)$ into the two terms $\frac{1}{2}g_i(s) f_i(x(s))$ and $\frac{1}{2}g_i(s) f_i(x(s))+e_i(t)$.
The first term is always nonnegative and used to carry out the estimates of the operator in the lemmas and corollaries in this section. We will make the second term
$\frac{1}{2}g_i(s) f_i(x(s))+e_i(t)$ nonnegative by choosing appropriate domains of $f_i$. The choice of the even split
of $g_i(s) f_i(x(s))$ here is not necessarily optimal in terms of obtaining maximal $\l$-intervals for
the existence of periodic solutions of the systems.

Noting that $g_i(t)$ is positive on $[0,T]$, $\lim_{ x \to 0} f_i(x) = \infty, i=1,...,n, $ implies that there is a $\delta>0$ such that
$$
f_i(x)\geq 2 \frac{\max_{t \in [0,T]} \{\abs{e_i(t)}+1\}}{\min_{t \in [0,T]}\{g_i(t)\}}, i=1,...,n,
$$
for $ x \in  \mathbb{R}_+^n, 0< \abs{x} \leq \delta$. Now for $x \in  \bar{\O}_r \setminus \{0\}$ and $0<r <\delta,$  noting that
$$
\delta > r \geq \sum_{i=1}^n |x_i(t)| \geq \min_{ t \in [0,T]} \sum_{i=1}^n |x_i(t)| \geq \s \norm{x}>0, \;\; t \in [0, T],
$$
and therefore, we have, for $t \in [0, T]$,
\begin{equation*}
\begin{split}
g_i(t)f_i(x(t))+e_i(t) &\geq \frac{1}{2}g_i(t)f_i(x(t))+e_i(t)\\
                       &\geq \frac{2}{2}g_i(t)\frac{\max_{t \in [0,T]} \{\abs{e_i(t)}+1\}}{\min_{t \in [0,T]}\{g_i(t)\}}+e_i(t)\\
                       & >0.
\end{split}
\end{equation*}
Thus, it is clear that $\mathcal{T}_{\l}^ix(t)$ in (\ref{T_def}) is well defined and positive, and now it is easy to see that  $\mathcal{T} _{\l}(\bar{\O}_r \setminus \{0\}) \subset K$
 and $\mathcal{T}_{\l}:\bar{\O}_r \setminus \{0\} \to K$ is completely continuous.

On the other hand, if $\lim_{ x \to \infty}  f_i(x)=\infty, i=1,\dots,n$, there is a $R''>0$ such that
$$
f_i(x)\geq 2 \frac{\max_{t \in [0,T]} \{\abs{e_i(t)}+1\}}{\min_{t \in [0,T]}\{g_i(t)\}}, i=1,...,n,
$$
for $ x \in  \mathbb{R}_+^n, \abs{x} \geq R''$. Now let $\Delta=\frac{R''}{\s}$. Then for $x \in  K\setminus \O_R,$ $R>\Delta$,  we have that $\min\limits_{0 \leq t \leq  T }\sum_{i=1}^nx_i(t) \geq \s \norm{x}\geq R'',$
and therefore,
$$g_i(t)f_i(x(t))+e_i(t)\geq \frac{1}{2}g_i(t)f_i(x(t))+e_i(t)>0, t \in [0,T].$$
Now $\mathcal{T}_{\l}^ix(t)$ in (\ref{T_def}) is well defined and positive. It is clear that  $\mathcal{T} _{\l}( K\setminus \O_R) \subset K$ and $\mathcal{T}_{\l}:K\setminus \O_R \to K$ is completely continuous.
\epf

Now let
$$
\G=\min_{i=1,...,n}\{ \frac{1}{2} m_i \s \int^{ T }_{0} g_i(s)ds\}>0.
$$
\slm\mlabel{lm5}
Assume \rm{(A),(H1),(H2)} hold and $e_i(t)\geq 0, t \in [0,T], i=1,...,n$. Let $r>0$ and if there exist $\eta>0$
and integer $j$, $1 \leq j \leq n$ such that
$$f_j(x(t)) \geq \eta \sum_{i=1}^n x_i(t) \for  t \in [0, T ],$$
for $ x(t)=(x_1(t),...,x_n(t)) \in \partial \O_r $, then the following inequality holds,
$$
\norm{\mathcal{T}_{\l}x} \geq \l \G \eta \norm{x}.
$$
\elm
\pf From the definition of $\mathcal{T}_{\l}x$ it follows that
\begin{equation*}
\begin{split}
 \norm{\mathcal{T}_{\l}x} & \geq\max\limits_{0 \leq t \leq T }\mathcal{T}_{\l}^j x(t)\\
 & \geq   \frac{1}{2} \l m_j\int^{ T }_{0} g_j(s) f_j(x(s))ds\\
 & \geq   \frac{1}{2} \l m_j \int^{ T }_{0} g_j(s)\eta \sum_{i=1}^n x_i(s)ds\\
 & \geq    \l m_j \frac{1}{2} \s \int^{ T }_{0} g_j(s)ds \eta \norm{x}\\
 & =  \l \G \eta \norm{x}.
\end{split}
\end{equation*}
\epf
If $e_i$ takes negative values, we need to adjust $\delta$ and $\Delta $ in Corollary \ref{cordefinition} to guarantee that $g_i(s)f_i(x(s))+ e_i(s)$ is nonnegative.
\scor\mlabel{lm51}
Assume \rm{(A),(H1),(H2), (H3)} hold. \\
(a). If $\lim_{ x \to 0}  f_i(x)=\infty, i=1,\dots,n,$  then Lemma \ref{lm5} is true if, in addition, $0<r< \delta$, where $\delta$ is defined in Corollary \ref{cordefinition}.\\
(b). If $\lim_{ \abs{x} \to \infty}  f_i(x)=\infty, i=1,\dots,n$,  then Lemma \ref{lm5} is true if, in addition, $ r> \Delta$,  where $\Delta$ is defined in Corollary \ref{cordefinition}.
\ecor
\pf
We split $g_i(s) f_i(x(s))+e_i(t)$ into the two terms $\frac{1}{2}g_i(s) f_i(x(s))$ and $\frac{1}{2}g_i(s) f_i(x(s))+e_i(t)$.
By choosing $\delta$ and $\Delta $ in Corollary \ref{cordefinition}, $g_i(s) f_i(x(s))+e_i(t)$ becomes
nonnegative. The estimate in Corollary \ref{lm51} can be carried out by the first term as in Lemma \ref{lm5}.
\epf

Let $\hat{f}_i(\theta): [1, \infty) \to \mathbb{R}_+$ be the function given by
$$\hat{f}_i(\theta) =\max \{f_i(u):u \in \mathbb{R}_+^n \; \rm{and}\; 1 \leq \abs{u} \leq \theta \}, i=1,...,n.$$
It is easy to see that $\hat{f}_i(\theta)$ is a nondecreasing function on $[1,\infty)$. The following lemma is essentially the same  as Lemma 2.8 in \cite{HWJMAA1}.
The following proof is only for completeness.
\slm (\cite{HWJMAA1})\mlabel{lm6}
Assume \rm{(H1)} holds. If $\lim_{ \abs{x} \to \infty} \frac{f_i(x)}{\abs{x}}$ exists (which can be infinity), then $\lim_{\theta \to \infty} \frac{\hat{f}_i(\theta)}{\theta}$ exists and
$\lim_{\theta \to \infty} \frac{\hat{f}_i(\theta)}{\theta}=\lim_{ \abs{x} \to \infty} \frac{f_i(x)}{\abs{x}}$.
\elm
\pf
We consider the two cases, (a) $f_i(x)$ is bounded for $\abs{x} \geq 1$ and (b) $f_i(x)$ is unbounded for $\abs{x} \geq 1$. For case (a),
it follows that $\lim_{\theta \to \infty} \frac{\hat{f}_i(\theta)}{\theta}=\lim_{ \abs{x} \to \infty} \frac{f_i(x)}{\abs{x}}=0$.
\noindent For  case (b), for any $\delta > 1$, let $M^i=\hat{f}_i(\delta)$
and
$$ N_{\delta}^i = \inf\{\abs{x}: x \in \mathbb{R}_+^n, \;\abs{x} \geq \delta, f_i(x) \geq M^i \} \geq \delta>1,
$$
then
$$
\max\{f_i(x): 1 \leq \abs{x} \leq N_{\delta}^i, \; x \in \mathbb{R}_+^n \}=M^i=\max\{f_i(x): \abs{x} = N_{\delta}^i, \; x \in \mathbb{R}_+^n \}.
$$
Therefore, for any $\delta > 1$, there exists a $ N_{\delta}^i \geq \delta$ such that
$$ \hat{f}_i(\theta) = \max\{f_i(x): N_{\delta}^i \leq \abs{x} \leq \theta , \; x \in \mathbb{R}_+^n \} \;\; {\rm for} \;\; \theta > N_{\delta}^i.
$$
Now, suppose that $b_i=\lim_{ \abs{x} \to \infty} \frac{f_i(x)}{\abs{x}}< \infty$. In other words, for any $\e>0$, there is a $\delta>1$ such that
\begin{equation}\label{e-d-defintion}
b_i -\e <\frac{f_i(x)}{\abs{x}} < b_i +\e,  \textrm{  for }  x \in \mathbb{R}_+^n,\; \abs{x} > \delta.
\end{equation}
Thus, for $\theta>N_{\delta}^i$, there exist $x_1, x_2 \in \mathbb{R}_+^n$ such that  $\abs{x_1}=\theta$, $\theta \geq \abs{x_2} \geq N_{\delta}^i$  and $f_i(x_2)=\hat{f}_i(\theta)$. Therefore,
\begin{equation}\label{e-d-inq}
\frac{f_i(x_1)}{\abs{x_1}} \leq  \frac{\hat{f}_i(\theta)}{\theta} = \frac{f_i(x_2)}{\theta} \leq \frac{f_i(x_2)}{\abs{x_2}}.
\end{equation}
(\ref{e-d-defintion}) and (\ref{e-d-inq}) yield that
\begin{equation}\label{e-d-new}
b_i -\e <\frac{\hat{f}_i(\theta)}{\theta} < b_i +\e \textrm{ for }  \theta>N_{\delta}^i.
\end{equation}
Hence $\lim_{\theta \to \infty} \frac{\hat{f}_i(\theta)}{\theta}=\lim_{ x \to \infty} \frac{f_i(x)}{\abs{x}}$. Similarly, we can show $\lim_{\theta \to \infty} \frac{\hat{f}^i(\theta)}{\theta}=\lim_{ x \to \infty} \frac{f_i(x)}{\abs{x}}$ if $\lim_{ \abs{x} \to \infty} \frac{f_i(x)}{\abs{x}}=\infty$.
\epf

\slm\mlabel{lm7}
Assume \rm{(A),(H1),(H2)} hold and $e_i(t)\geq 0, t \in [0,T], i=1,...,n$.\\ Let $r > \max \{\frac{1}{\s}, 2 \l \sum_{i=1}^{n } M_i\int_0^T \abs{e_i(s)}ds \}$ and if there exists an $\e > 0$ such that
$$
\hat{f}_i(r) \leq  \e r,\;\; i=1,...,n,
$$
then
$$
\norm{\mathcal{T}_{\l}x} \leq  \l \hat{C}\e \norm{x} + \frac{1}{2} \norm{x} \;\; {\rm for} \;\; x \in \partial\O_{r}.
$$
\elm
where the constant $\hat{C}=\sum_{i=1}^{n } M_i \int^{ T}_0 g_i(s)ds$.
\pf From the definition of $\mathcal{T}_{\l}$, we have for $x \in \partial\O_{r}$,
 \begin{eqnarray*}
  \norm{\mathcal{T}_{\l}x} & = & \sum_{i=1}^{n }\max\limits_{0 \leq t \leq T }\mathcal{T}_{\l}^i x(t)  \\
 & \leq &  \sum_{i=1}^{n }\l M_i \int^{ T}_0 g_i(s)f_i(x(s))ds + \l M_i \sum_{i=1}^{n } \int_0^T \abs{e_i(s)}ds \\[.2cm]
  & \leq &   \sum_{i=1}^{n } \l M_i \int^{ T}_0 g_i(s)\hat{f}_i(r)ds +  \frac{ r}{2} \\[.2cm]
  & \leq &  \sum_{i=1}^{n } \l M_i \int^{ T}_0 g_i(s)ds r\e +   \frac{r}{2}\\
  & = & \l \hat{C}\e \norm{x} + \frac{1}{2} \norm{x}.
\end{eqnarray*}
\epf

If $e_i$ takes negative values, we need to restrict the domain of $\mathcal{T}_{\l}$ to guarantee that $g_i(s)f_i(x(s))+ e_i(s)$ is nonnegative.
\scor\mlabel{lm71}
Assume \rm{(A),(H1),(H2), (H3)} hold. If $\lim_{ x \to \infty}  f_i(x)=\infty, i=1,\dots,n$, Lemma \ref{lm7} is true if , in addition, $r> \Delta$,  where $\Delta$ is defined in Corollary \ref{cordefinition}.
\ecor
\pf
If we choose $\Delta$ defined in Corollary \ref{cordefinition}, then $\mathcal{T}_{\l}$ is well-defined and $g_i(s)f_i(x(s))+ e_i(s)$ is nonnegative,
and Corollary \ref{lm71} can be shown
in the same way as Lemma \ref{lm7}.
\epf
The conclusions of Lemmas \ref{lm5} and \ref{lm7} are based on the inequality assumptions between $f(x)$ and $x$. If
these assumptions is not necessarily true, we will have the following results.
\slm\mlabel{lm8}
Assume \rm{(A),(H1),(H2)} hold and $e_i(t)\geq 0, t \in [0,T], i=1,...,n$. Let $r>0$. Then
$$
\norm{\mathcal{T}_{\l}x}  \geq \l \sum_{i=1}^{n } \frac{m_i \hat{m}_{r}}{2}   \int^{ T }_{0} g_i(s)ds,
$$
for all $x \in \partial \O_{r}$,
where $\hat{m}_{r}=\min \{ f_i(x):  x \in \mathbb{R}_+^n \; \rm{and}\; \s r \leq \abs{x} \leq r, \,\, i=1,...,n\}>0.$
\elm
\pf  If  $ x(t) \in \partial \O_{r}$, then $\s r \leq \abs{x(t)}= \sum_{i=1}^n |x_i(t)| \leq r, t \in [0, T]$.  Therefore $ f_{i}(x(t)) \geq  \hat{m}_{r} \; \rm{for}\; t \in [0, T ],$
$i=1,...,n$. By the definition of $\mathcal{T}_{\l}$, we have
\begin{equation*}
\begin{split}
 \norm{\mathcal{T}_{\l}x} & = \sum_{i=1}^{n } \max\limits_{0 \leq t \leq T }\mathcal{T}_{\l}^i x(t)\\
 & \geq   \sum_{i=1}^{n } \frac{1}{2} \l m_i\int^{ T }_{0} g_i(s) f_i(x(s))ds\\
 & \geq  \l \sum_{i=1}^{n } \frac{m_i\hat{m}_{r}}{2}   \int^{ T }_{0} g_i(s) ds.\\
\end{split}
\end{equation*}
\epf
Now we consider the case that $e_i$ may take negative values.  We need to restrict the domain of $\mathcal{T}_{\l}$ to guarantee that $g_i(s)f_i(x(s))+ e_i(s)$ is nonnegative.
$\frac{1}{2}g_i(s) f_i(x(s))$ is used to carry out the estimates in Lemma \ref{lm8}.

\scor\mlabel{lm81}
Assume \rm{(A),(H1),(H2), (H3)} hold.\\
(a). If $\lim_{ x \to 0}  f_i(x)=\infty, i=1,\dots,n$, Lemma \ref{lm8} is true if, in addition,  $0<r <\delta$, where $\delta>0$ is defined in Corollary \ref{cordefinition}.\\
(b). If $\lim_{ \abs{x} \to \infty}  f_i(x)=\infty, i=1,\dots,n$, Lemma \ref{lm8} is true if, in addition, $r  > \Delta$, where $\Delta>0$ is defined in Corollary \ref{cordefinition}.
\ecor
\pf
By selecting  $\delta$ and $\Delta$ defined in Corollary \ref{cordefinition}, $\mathcal{T}_{\l}$ is well-defined and $g_i(s)f_i(x(s))+ e_i(s)$ is nonnegative, and then
Corollary \ref{lm81} can be shown
as Lemma \ref{lm8}.
\epf

\slm\mlabel{lm9}
Assume \rm{(A),(H1),(H2)} hold and $e_i(t)\geq 0, t \in [0,T], i=1,...,n$. Let $r>0$. Then
$$
\norm{\mathcal{T}_{\l}x} \leq   \l \big( \sum_{i=1}^{n } M_i \int^{ T}_0 g_i(s)\hat{M}_{r} ds +  \sum_{i=1}^{n } M_i\int_0^T \abs{e_i(s)}ds \big ),
$$
for all $ x \in \partial \O_{r}$,
where $\hat{M}_{r}=\max\{ f_i(u):  u \in \mathbb{R}_+^n \; \rm{and}\; \s r\leq \abs{u} \leq r, \,\, i=1,...,n\}>0$.
\elm
\pf If  $ x \in \partial \O_{r}$, then $\s r \leq \abs{x(t)} \leq r, t \in [0, T]$.  Therefore $ f_{i}(x(t)) \leq \hat{M}_{r} \; \rm{for}\; t \in [0, T]$, $i=1,...,n.$  Thus we have that
 \begin{eqnarray*}
  \norm{\mathcal{T}_{\l}x} & = & \sum_{i=1}^{n }\max\limits_{0 \leq t \leq T }\mathcal{T}_{\l}^i x(t)  \\
 & \leq &  \sum_{i=1}^{n }\l M_i \int^{ T}_0 g_i(s)f_i(x(s))ds + \l \sum_{i=1}^{n } M_i\int_0^T \abs{e_i(s)}ds \\[.2cm]
  & \leq &  \sum_{i=1}^{n }\l M_i \int^{ T}_0 g_i(s)f_i(x(s))ds + \l \sum_{i=1}^{n } M_i\int_0^T \abs{e_i(s)}ds \\[.2cm]
  & \leq &   \sum_{i=1}^{n } \l M_i \int^{ T}_0 g_i(s)\hat{M}_{r} ds +  \l \sum_{i=1}^{n } M_i\int_0^T \abs{e_i(s)}ds  \\[.2cm]
  & \leq &  \l \big( \sum_{i=1}^{n } M_i \int^{ T}_0 g_i(s)\hat{M}_{r} ds +  \sum_{i=1}^{n } M_i\int_0^T \abs{e_i(s)}ds \big ).
\end{eqnarray*}
\epf
Again, if $e_i$ takes negative values, we need to restrict $r$ and $R$ to guarantee $g_i(s)f_i(x(s))+ e_i(s)$ is nonnegative.
\scor\mlabel{lm91}
Assume \rm{(A),(H1),(H2), (H3)} hold.\\
(a). If $\lim_{ x \to 0}  f_i(x)=\infty, i=1,\dots,n$, Lemma \ref{lm9} is true if, in addition,  $0<r < \delta$, where $\delta>0$ is defined in Corollary \ref{cordefinition}.\\
(b). If $\lim_{ x \to \infty}  f_i(x)=\infty, i=1,\dots,n$, Lemma \ref{lm9} is true if, in addition,  $r> \Delta$, where $\Delta>0$ is defined in Corollary \ref{cordefinition}.
\ecor
\pf
By selecting  $\delta$ and $\Delta$ defined in Corollary \ref{cordefinition}, $\mathcal{T}_{\l}$ is well-defined and $g_i(s)f_i(x(s))+ e_i(s)$ is nonnegative, and then the Corollary can be shown exactly
as Lemma \ref{lm9}.
\epf

\section{Proof of Theorem \ref{th3}}

\pf Part (a). Since $e_i(t) \geq 0$, $\mathcal{T}_{\l}$ is defined on $K\setminus \{0\}$ and $g_i(s)f_i(x(s))+ e_i(s)$ is
nonnegative.
Noting $\lim_{ \abs{x} \to \infty} \frac{f_i(x)}{\abs{x}}=0$,  $i=1,\dots,n$,  it follows from Lemma ~\ref{lm6} that
$\lim_{\theta \to \infty} \frac{\hat{f}_i(\theta)}{\theta}=0$, $i=1,...,n.$ Therefore, we can choose $r_1 > \max \{\frac{1}{\s}, 2 \l \sum_{i=1}^{n } M_i\int_0^T \abs{e_i(s)}ds \}$
so that $\hat{f}_i(r_1) \le \e  r_1 ,$ $i=1,...,n$, where the constant $\e> 0$ satisfies
$$
\l \hat{C}\e  < \frac{1}{2},
$$
and $\hat{C}$ is the positive constant defined in Lemma \ref{lm7}. We have by Lemma ~\ref{lm7} that
$$
\norm{\mathcal{T}_{\l}x} \leq (\l \hat{C}\e + \frac{1}{2}) \norm{x} < \norm{x} \quad \textrm{for} \quad  x \in \partial\O_{r_1}.
$$
On the other hand, by the condition $\lim_{ x \to 0} f_i(x)=\infty$, there there is a positive number $r_2 <  r_1$ such that
$$
f_i(x) \geq \eta \abs{x}, \;\;i=1,...,n
$$
for $ x=(x_1,...,x_n) \in \mathbb{R}_+^n\setminus \{0\} $ and $ \abs{x} \leq r_2,$
where $\eta > 0$ is chosen so that
$$
\l \G \eta > 1.
$$
It is easy to see that, for $x=(x_1,...,x_n) \in  \partial \O_{r_2}, \;\;t \in [0,T],$
$$f_i(x(t)) \geq \eta\sum_{i=1}^n x_i(t).$$
Lemma ~\ref{lm5} implies that
$$
\norm{\mathcal{T}_{\l}x} \geq \l \G \eta \norm{x} > \norm{x} \quad \textrm{for} \quad  x \in \partial\O_{r_2}.
$$
By Lemma ~\ref{lm1}, $\mathcal{T}_{\l}$ has a fixed point $x \in  \bar{\O}_{r_1} \setminus \O_{r_2}$.
The fixed point $x \in \bar{\O}_{r_1} \setminus \O_{r_2}$ is the desired positive periodic solution of (\ref{eq3}).

Part (b). Again since $e_i(t) \geq 0$, $\mathcal{T}_{\l}$ is defined on $K \setminus \{0\}$ and $g_i(s)f_i(x(s))+ e_i(s)$ is
nonnegative. Fix two numbers $0 < r_{3} < r_{4},$
there exists a $\l_0 >0$ such that
$$
\l_0  < \frac{r_3}{ \sum_{i=1}^{n } M_i \int^{ T}_0 g_i(s)\hat{M}_{r_3} ds +  \sum_{i=1}^{n } M_i\int_0^T \abs{e_i(s)}ds },
$$
and
$$
\l_0  < \frac{r_4}{ \sum_{i=1}^{n } M_i \int^{ T}_0 g_i(s)\hat{M}_{r_4} ds +  \sum_{i=1}^{n } M_i\int_0^T \abs{e_i(s)}ds },
$$
where $\hat{M}_{r_3}$ and $\hat{M}_{r_4}$ are defined in Lemma ~\ref{lm9}.
Thus, Lemma \ref{lm9} implies that, for $ 0< \l < \l_{0} $,
$$
\norm{\mathcal{T}_{\l} x} < \norm{x}, \;\;\rm{for} \;\; x \in  \partial
\O_{r_j}, \;\; (j=3,4).
$$
On the other hand, in view of the assumptions $\lim_{ x \to \infty} \frac{f_i(x)}{\abs{x}}=\infty$ and $\lim_{ x \to 0} f_i(x)=\infty$, there are positive numbers $0< r_2 < r_3< r_4< r'_1$ such that
$$
f_i(x) \geq \eta \abs{x}
$$
for $ x=(x_1,...,x_n) \in \mathbb{R}_+^n$ and $ 0 <\abs{x} \leq  r_2$ or $ \abs{x} \geq  r'_1$
where $\eta > 0$ is chosen so that
$$
\l \G \eta > 1.
$$
Thus if  $ x=(x_1,...,x_n) \in \partial \O_{r_2}$, then
$$f_i(x(t)) \geq \eta\sum_{i=1}^n x_i(t), \;\; t \in [0, T].$$
Let $r_1 = \max\{2r_4, \frac{1}{\s} r'_1 \}$. If $ x=(x_1,...,x_n) \in \partial \O_{r_1}$, then
$$ \min_{ 0 \leq t \leq  T} \sum_{i=1}^n x_i(t) \geq \s
\norm{x}= \s r_1 \geq r'_1,$$
which implies that
$$f_i(x(t)) \geq \eta \sum_{i=1}^n x_i(t) \; \rm{for} \; t \in [0, T].
$$
Thus Lemma ~\ref{lm5} implies that
$$
\norm{\mathcal{T}_{\l}x} \geq \l \G \eta \norm{x} > \norm{x} \quad \textrm{for} \quad  x \in \partial\O_{r_1},
$$
and
$$
\norm{\mathcal{T}_{\l}x} \geq \l \G \eta \norm{x} > \norm{x} \quad \textrm{for} \quad  x \in \partial\O_{r_2}.
$$
It follows from Lemma ~\ref{lm1}, $\mathcal{T}_{\l}$ has two fixed points $x_1(t)$ and $x_2(t)$ such that  $x_1(t) \in \bar{\O}_{r_3} \setminus \O_{r_2}$
and $x_2(t) \in \bar{\O}_{r_1} \setminus \O_{r_4}$ , which are
the desired distinct positive periodic solutions of (\ref{eq3}) for $ \l < \l_0$ satisfying
$$
r_1 < \norm{x_1} < r_3 < r_4 < \norm{x_2} < r_2.
$$

Part (c). First we note that $\mathcal{T}_{\l}$ is defined on $K \setminus \{0\}$ and $g_i(s)f_i(x(s))+ e_i(s)$ is
nonnegative since $e_i(t) \geq 0$.  Fix a number $r_{3}>0 $. Lemma \ref{lm9} implies that there exists a $\l_0 >0$ such that we have, for $ 0< \l < \l_{0} $,
$$
\norm{\mathcal{T}_{\l}x} < \norm{x}, \;\;\rm{for} \;\; x \in  \partial
\O_{r_3}.
$$
On the other hand,  in view of the assumption $\lim_{ x \to 0} f_i(x)=\infty$, there is a positive number $0< r_2 < r_3$ such that
$$
f_i(x) \geq \eta \abs{x}
$$
for $ x=(x_1,...,x_n) \in \mathbb{R}_+^n$ and $ 0< \abs{x} \leq  r_2$
where $\eta > 0$ is chosen so that
$$
\l \G \eta > 1.
$$
Thus if  $ x=(x_1,...,x_n) \in \partial \O_{r_2}$, then
$$f_i(x(t)) \geq \eta\sum_{i=1}^n x_i(t), \;\; t \in [0, T].$$
Thus Lemma ~\ref{lm5} implies that
$$
\norm{\mathcal{T}_{\l}x} \geq \l \G \eta \norm{x} > \norm{x} \quad \textrm{for} \quad  x \in \partial\O_{r_2}.
$$
Lemma ~\ref{lm1} implies that $\mathcal{T}_{\l}$ has a fixed point $x \in  \bar{\O}_{r_3} \setminus \O_{r_2}$.
The fixed point $x \in \bar{\O}_{r_3} \setminus \O_{r_2}$ is the desired positive periodic solution of (\ref{eq3}).
\epf

\section{Proof of Theorem \ref{th5}}
%

\pf
Part (a). Since $\lim_{ \abs{x} \to \infty}  f_i(x)=\infty, i=1,\dots,n$, by Corollary \ref{cordefinition} , there is a $\Delta>0$ such that if $R>\Delta$,
then $g_i(s)f_i(x(s))+ e_i(s)$ is nonnegative and $\mathcal{T}_{\l}: K\setminus \O_R \to K$ is defined.
Now for a fixed number $r_{1}>\Delta$,  Corollary \ref{lm81} implies that there exists a $\l_0 >0$
such that, for $ \l > \l_{0} $,
$$
\norm{\mathcal{T}_{\l}x} > \norm{x}, \;\;\rm{for} \;\; x \in  \partial
\O_{r_1}.
$$
On the other hand, since $\lim_{ \abs{x} \to \infty} \frac{f_i(x)}{\abs{x}}=0$,  $i=1,\dots,n$,  it follows from Lemma ~\ref{lm6} that
$\lim_{\theta \to \infty} \frac{\hat{f}_i(\theta)}{\theta}=0$, $i=1,...,n.$ Therefore, we can choose $$r_2 > \max \{2r_1, \frac{1}{\s}, 2 \l \sum_{i=1}^{n } M_i\int_0^T \abs{e_i(s)}ds \}>\Delta$$
so that $\hat{f}_i(r_2) \le \e  r_2 ,$ $i=1,...,n$, where the constant $\e> 0$ satisfies
$$
\l \hat{C}\e  < \frac{1}{2},
$$
We have, by Corollary ~\ref{lm71},  that
$$
\norm{\mathcal{T}_{\l}x} \leq (\l \hat{C}\e + \frac{1}{2}) \norm{x} < \norm{x} \quad \textrm{for} \quad  x \in \partial\O_{r_2}.
$$
By Lemma ~\ref{lm1}, $\mathcal{T}_{\l}$ has a fixed point $x \in  \bar{\O}_{r_2} \setminus \O_{r_1}$.
The fixed point $x \in \bar{\O}_{r_2} \setminus \O_{r_1}$ is the desired positive periodic solution of (\ref{eq3}).

Part (b).  First, since $\lim_{ x \to 0}  f_i(x)=\infty, i=1,...,n$,
by Corollary \ref{cordefinition}, there is $\delta>0$ such that if $0< r< \delta$, $\mathcal{T}_{\l}$
is defined on $\bar{\O}_r \setminus \{0\}$ and $g_i(s)f_i(x(s))+ e_i(s)$ is nonnegative.
Furthermore, $\mathcal{T} _{\l}(\bar{\O}_r \setminus \{0\}) \subset K$.
Now for a fixed number $r_{1}<\delta$,  and Corollary \ref{lm91} implies that there exists a $\l_1 >0$
such that we have, for $ \l < \l_{1} $,
$$
\norm{\mathcal{T}_{\l}x} < \norm{x}, \;\;\rm{for} \;\; x \in  \partial
\O_{r_1}.
$$
In view of the assumption $\lim_{ x \to 0} f_i(x)=\infty$, there is a positive number $0< r_3 < r_1$ such that
$$
f_i(x) \geq \eta \abs{x}
$$
for $ x=(x_1,...,x_n) \in \mathbb{R}_+^n$ and $ 0< \abs{x} \leq  r_3$
where $\eta > 0$ is chosen so that
$$
\l \G \eta > 1.
$$
Thus if  $ x=(x_1,...,x_n) \in \partial \O_{r_3}$, then
$$f_i(x(t)) \geq \eta\sum_{i=1}^n x_i(t), \;\; t \in [0, T].$$
Thus Corollary ~\ref{lm51} implies that
$$
\norm{\mathcal{T}_{\l}x} \geq \l \G \eta \norm{x} > \norm{x} \quad \textrm{for} \quad  x \in \partial\O_{r_3}.
$$
It follows from Lemma ~\ref{lm1}, $\mathcal{T}_{\l}$ has a fixed point $x_1(t)$ such that  $x_1(t) \in \bar{\O}_{r_1} \setminus \O_{r_3}$
which is a positive periodic solutions of (\ref{eq3}) for $ \l < \l_1$ satisfying
$$
r_3 < \norm{x_1} < r_1.
$$
On the other hand, Since $\lim_{ \abs{x} \to \infty} \frac{f_i(x)}{\abs{x}}=\infty, i=1,\dots,n$,
by Corollary \ref{cordefinition}, there is $\Delta>0$ such that if $R>\Delta$, $\mathcal{T}_{\l}$
is defined on $K\setminus \O_R$ and $g_i(s)f_i(x(s))+ e_i(s)$ is nonnegative. Furthermore, $\mathcal{T} _{\l}(K  \setminus \O_R) \subset K$.
For a fixed number $r_{2}> \max\{\Delta, r_1\}$,  and Corollary \ref{lm91} implies that there exists a $0< \l_0 < \l_1$
such that we have, for $ \l < \l_{0} $,
$$
\norm{\mathcal{T}_{\l}x} < \norm{x}, \;\;\rm{for} \;\; x \in  \partial
\O_{r_2}.
$$
Since $\lim_{ \abs{x} \to \infty} \frac{f_i(x)}{\abs{x}}=\infty, \;\; i=1,\dots,n, $ implies that there is a positive number $r'$
such that
$$
f_i(x) \geq \eta \abs{x}
$$
for $ x=(x_1,...,x_n) \in \mathbb{R}_+^n$ and $ \abs{x} \geq  r'$
where $\eta > 0$ is chosen so that
$$
\l \G \eta > 1.
$$
Let $r_4 = \max\{2r_2, \frac{1}{\s} r' \}>\Delta$. If $ x=(x_1,...,x_n) \in \partial \O_{r_4}$, then
$$ \min_{ 0 \leq t \leq  T} \sum_{i=1}^n x_i(t) \geq \s
\norm{x}= \s r_4 \geq r',$$
which implies that
$$f_i(x(t)) \geq \eta \sum_{i=1}^n x_i(t) \; \rm{for} \; t \in [0, T].
$$
Again Corollary ~\ref{lm51} implies that
$$
\norm{\mathcal{T}_{\l}x} \geq \l \G \eta \norm{x} > \norm{x} \quad \textrm{for} \quad  x \in \partial\O_{r_4}.
$$
It follows from Lemma ~\ref{lm1}, $\mathcal{T}_{\l}$ has a fixed point $x_2(t) \in \bar{\O}_{r_4} \setminus \O_{r_2}$ , which is a positive periodic
solutions of (\ref{eq3}) for $ \l < \l_0$ satisfying
$$
r_2 < \norm{x_2} < r_4.
$$
Noting that
$$
r_3 < \norm{x_1} < r_1 <r_2< \norm{x_2} < r_4,
$$
we can conclude that $x_1$ and $x_2$ are the desired distinct positive periodic solutions of (\ref{eq3}) for $ \l < \l_0$.

Part (c).  Since $\lim_{ x \to 0}  f_i(x)=\infty, i=1,\dots,n$, by Corollary \ref{cordefinition} , there is a $\delta>0$ such that if $0<r <\delta$,
then $\mathcal{T}_{\l}$ is defined and $g_i(s)f_i(x(s))+ e_i(s)$ is nonnegative.
Now for a fixed number $r_{1}< \delta $,  Corollary \ref{lm91} implies that there exists a $\l_1 >0$
such that we have, for $ \l < \l_{1} $,
$$
\norm{\mathcal{T}_{\l}x} < \norm{x}, \;\;\rm{for} \;\; x \in  \partial
\O_{r_1}.
$$
On the other hand,  in view of the assumption $\lim_{ x \to 0} f_i(x)=\infty$, there there is a positive number $0< r_2 < r_1<\delta$ such that
$$
f_i(x) \geq \eta \abs{x}
$$
for $ x=(x_1,...,x_n) \in \mathbb{R}_+^n$ and $ 0< \abs{x} \leq  r_2$
where $\eta > 0$ is chosen so that
$$
\l \G \eta > 1.
$$
Thus if  $ x=(x_1,...,x_n) \in \partial \O_{r_2}$, then
$$f_i(x(t)) \geq \eta\sum_{i=1}^n x_i(t), \;\; t \in [0, T].$$
Thus Corollary ~\ref{lm51} implies that
$$
\norm{\mathcal{T}_{\l}x} \geq \l \G \eta \norm{x} > \norm{x} \quad \textrm{for} \quad  x \in \partial\O_{r_2}.
$$
Lemma ~\ref{lm1} implies that $\mathcal{T}_{\l}$ has a fixed point $x_1 \in  \bar{\O}_{r_1} \setminus \O_{r_2}$.
The fixed point $x_1 \in \bar{\O}_{r_1} \setminus \O_{r_2}$ is the desired positive periodic solution of (\ref{eq3}).
\epf

\section{Acknowledgment}
The author would like to thank Jifeng Chu
for bringing several references to his attention.

\end{document}